\documentclass[12pt]{article}
\usepackage{amsmath,amsthm,amssymb}
\usepackage{fullpage}
\usepackage{enumerate}

\renewcommand{\pmod}[1]{\,(\textup{mod}\,#1)}
\numberwithin{equation}{section}
\theoremstyle{plain}
\newtheorem{theorem}{Theorem}[section]

\newtheorem{lemma}[theorem]{Lemma}
\newtheorem{corollary}[theorem]{Corollary}

\newcommand{\mP}{\mathcal{P}}
\newcommand{\mPt}{\widetilde{\mathcal{P}}}
\newcommand{\mEt}{\widetilde{\mathcal{E}}}

\newcommand{\mQ}{\mathcal{Q}}
\newcommand{\mE}{\mathcal{E}}
\newcommand{\mD}{\mathcal{D}}
\newcommand{\SL}{\text{SL}_2(\mathbb{Z})}
\newcommand{\Sk}{S_k(\Gamma)}
\newcommand{\Mk}{M_k(\Gamma)}
\newcommand{\Minf}{M_k^{\infty}(\Gamma)}
\newcommand{\GN}{\Gamma}
\newcommand{\UH}{\mathcal{H}}

\newcommand{\D}{\Delta}

\begin{document}

\title{\Large \bf Integrable systems and modular forms of level 2}
\author{Mark J Ablowitz$^1$, Sarbarish Chakravarty$^2$ and 
Heekyoung Hahn$^3$ \\[1ex]
\small$^1$\it 
Department of Applied Mathematics, University of Colorado, Boulder, CO, 80309\\
\small$^2$\it
Department of Mathematics, University of Colorado, Colorado Springs, CO, 80933\\
\small$^3$\it
Department of Mathematics, University of Rochester, Rochester, NY, 14620}
\date{}

\maketitle

\begin{abstract}
\noindent A set of nonlinear differential equations associated with
the Eisenstein series of the congruent subgroup $\GN_0(2)$ of the
modular group $\SL$ is constructed. These nonlinear equations are 
analogues of the well known
Ramanujan equations, as well as the Chazy and 
Darboux-Halphen equations associated
with the modular group. The general solutions of these equations can be
realized in terms of the Schwarz triangle function $S(0,0,1/2; z)$.
\end{abstract}


\thispagestyle{empty}
\newpage

\section{Introduction}

In 1881, G. Halphen considered the nonlinear differential system \cite{Halphen}
\begin{align}
u_1'+u_2'&=u_1u_2, \nonumber\\
u_2'+u_3'&=u_2u_3,\label{halphen} \\
u_3'+u_1'&=u_3u_1 \,,\nonumber
\end{align}
for the functions $u_1(z),\, u_2(z), \, u_3(z)$, which originally appeared in 
Darboux's work of triply orthogonal surfaces on $\mathbb{R}^3$ \cite{Darboux}.
Note that $f'$ indicates derivation with respect to the argument
of the function $f$ throughout this article.
Halphen expressed the solution of the system \eqref{halphen} in terms of the logarithmic
derivatives of the null theta functions; namely,
$$
u_1(z) = \big(\ln \vartheta_4(0|z)\big)'\,, \qquad
u_2(z) = \big(\ln \vartheta_2(0|z)\big)'\,, \qquad
u_3(z) = \big(\ln \vartheta_3(0|z)\big)'\,,
$$
where the null theta functions are defined as (see e.g., \cite{WW, Rankin1})
\begin{equation*}
\vartheta_2(0|z):= \sum_{n=-\infty}^{\infty}q^{\frac{1}{2}(n+\frac{1}{2})^2}\,,
\qquad
\vartheta_3(0|z):= \sum_{n=-\infty}^{\infty}q^{\frac{1}{2}n^2}\,,
\qquad
\vartheta_4(0|z):= \sum_{n=-\infty}^{\infty}(-1)^nq^{\frac{1}{2}n^2}\,,
\end{equation*}
and $q:=e^{2 \pi iz}, \,\, \text{Im }z>0$.  
In 1909, Chazy in his study of Painlev\'e type equations of third order,
considered the nonlinear differential equation for the complex function $y(z)$
\begin{equation}\label{chazy}
y'''=2yy''-3y'\,^2,
\end{equation}
which as he noted, is related to the Darboux-Halphen (DH) system 
via $y=u_1+u_2+u_3$ \cite{Chazy}. 
It turns out that both the Chazy and the DH equations are intimately
connected to the theory of modular forms \cite{Ablowitz-T}. Indeed, a particular
solution of \eqref{chazy} is given by
\begin{equation*}
y(z):=\pi i E_2(z),
\end{equation*}
where
\begin{equation}\label{E}
E_2(z)=1-24\sum_{n=1}^{\infty}\Big(\sum_{d|n}d\Big)q^n \,,
\end{equation}
with $q:=e^{2 \pi iz}$, and $\text{Im }z>0$. 
$E_2(z)$ is the weight $2$ Eisenstein series associated with the
full modular group $\SL.$ The normalized Eisenstein series on $\SL$ are
defined, for even integer $k\geq 2,$ by
\begin{equation}
E_k(z)=1-\frac{2k}{B_k} \sum_{n=1}^{\infty}\sigma_{k-1}(n)q^n,
\label{Ek}
\end{equation}
where $B_k$ is the $k$th Bernoulli number and
\begin{equation*} 
\sigma_k(n):=\sum_{d|n}d^k \,, \quad n \in \mathbb{N}\,.
\end{equation*}

The purpose of the present paper is to find nonlinear ordinary
differential equations (ODEs) similar to the Halphen and Chazy equations,
and which are related to modular forms associated with subgroups of 
the modular group. In particular, we construct 
nonlinear ODEs satisfied by the Eisenstein series associated with the 
subgroup $\GN_0(2)$ (defined below) of the modular group $\SL$. We then illustrate
their relationships with certain system of ODEs also found by Halphen \cite{Halphen1}, 
as well as scalar ODEs considered by Schwarz \cite{Schwarz}, and 
later by Bureau \cite{Bureau}. The latter class of equations are known to be
important in the theory of conformal mapping \cite{Nehari}.
More recently such ODEs have appeared in several areas of mathematical physics
including magnetic monopoles \cite{Atiyah-H}, self-dual Yang-Mills and Einstein
equations \cite{Ablowitz-H3, Hitchin}, as well as topological field theory \cite{Dubrovin}.

Another motivation for this article stems from Ramanujan's work on modular forms, and 
subsequent extensions of some of Ramanujan's results by Ramamani.
In 1916, Ramanujan \cite{Ramanujan-arith}, \cite[pp 136--162]{Ramanujan-collect}, 
introduced the functions $P(q)$, $Q(q)$ and $R(q)$ defined for $|q|<1$ by
\begin{equation}
P(q):=1-24\sum_{n=1}^{\infty} \sigma_1(n)q^n\,, \quad
Q(q):=1+240\sum_{n=1}^{\infty} \sigma_3(n)q^n\,,\quad 
R(q):=1-504\sum_{n=1}^{\infty} \sigma_5(n)q^n \,, 
\label{PQR}
\end{equation}
and proved that these functions \eqref{PQR} satisfy the ODEs  
\begin{equation}
qP'=\frac{P^2-Q}{12}\,, \quad \qquad
qQ'=\frac{PQ-R}{3}\,, \quad \qquad
qR'=\frac{PR-Q^2}{2}\,.
\label{dPQR}
\end{equation}
We note here that the functions $P(q),\, Q(q)$ and $R(q)$ are the Eisenstein series
for $\SL$ introduced above. That is, $P(q)=E_2(z),$ $Q(q)=E_4(z),$ and $R(q)=E_6(z).$ 
There is an important correspondence between the Chazy equation and Ramanujan's
work. If we rewrite the system \eqref{dPQR} into a single equation, then it can 
be shown that $\pi i P(q)$ is a solution of the Chazy equation \cite{Ablowitz-H3}.

An analogue of the Ramanujan-type nonlinear system was
considered by Ramamani \cite{Ramamani-thesis}, who introduced three 
functions similar to \eqref{PQR}, and defined for $|q|<1$ by
\begin{equation}
\mP(q):=1-8\sum_{n=1}^{\infty}\frac{(-1)^{n}nq^n}{1-q^n}\,, \quad
\mPt(q):=1+24\sum_{n=1}^{\infty}\frac{nq^n}{1+q^n}\,, \quad 
\mQ(q):=1+16\sum_{n=1}^{\infty}\frac{(-1)^{n}n^3q^n}{1-q^n}\,. 
\label{mpptq}
\end{equation}
In the same manner as \eqref{dPQR}, these functions 
$\mP(q)$, $\mPt(q)$ and $\mQ(q)$ satisfy the differential equations 
\cite{Ramamani-thesis, Ramamani-paper, Hahn-eisen}
\begin{equation}
q\mP'=\frac{\mP^2-\mQ}{4}\,, \qquad
q\mPt'=\frac{\mPt\mP-\mQ}{2}\,, \qquad q\mQ'=\mP\mQ-\mPt\mQ \,. 
\label{dmpq}
\end{equation}
It turns out that the functions \eqref{mpptq} are the Eisenstein 
series associated with the congruence (congruent modulo 2) 
subgroup $\Gamma_0(2)$ of $\SL$.  Here $\GN_0(2)$ is defined by
$$\Gamma_{0}(2):=\left\{\gamma = \begin{pmatrix} a & b \\ c & d\end{pmatrix}
\in \SL : c \equiv 0 \, \pmod{2}\right\}.$$
In more contemporary notation, the normalized Eisenstein series on $\Gamma_0(2)$ 
are defined (see e.g., \cite{Schoeneberg}), for even integer $k\geq 2$, by
\begin{equation}\label{eisen}
\mE_k(z)=1+\frac{2k}{(1-2^k)B_k} \sum_{n=1}^{\infty}
\frac{(-1)^{n} n^{k-1}q^n}{1-q^n}\,,
\end{equation}
where, as before, $q:=e^{2 \pi iz}$ with $\text{Im }z>0$.
It is easy to check that $\mE_2(z)=\mP(q)$ and $\mE_4(z)=\mQ(q).$ If we 
define $\mEt_2(z)$ by
\begin{equation}\label{Et}
\mEt_2(z):=\mPt(q),
\end{equation}
then it can be shown that $2\mEt_2(z) = 3\mE_2(z)-E_2(z)$ 
(see Lemma \ref{3mp2mpt} in section \ref{sym}). Therefore, all three functions 
$\mP(q),\,\mPt(q)$, and $\mQ(q)$ introduced by Ramamani are expressible in terms of 
Eisenstein series. 

In this article, we show that Ramamani's system \eqref{dmpq} is equivalent to a third order
scalar nonlinear ODE found by Bureau in \cite{Bureau}, and whose solutions are given implicitly
by a Schwarz triangle function \cite{Schwarz, Nehari}. We also construct a Halphen-type ODE
system for $\GN_0(2)$ and show its equivalence to Ramamani's system. The solution to this
Halphen-type system is given in terms of the logarithmic derivatives involving the Schwarz
triangle function, in much the same way as the null theta functions solve \eqref{halphen}.
The paper is organized as follows.  In the following section, we introduce some basic 
definitions and preliminary materials on complex functions associated with the full modular 
group $\SL$, as well as its subgroup $\Gamma_0(2)$. In section 3, we discuss
the transformation properties of the modular forms on $\GN_0(2)$, and present our main results
in section 4. 

\section{Background}
The aim of this section is to provide
a brief background on modular forms necessary for the subsequent discussions 
of the various nonlinear ODEs on $\Gamma_0(2)$ and transformation properties 
of their solutions. Detailed discussions on modular forms
can be found in several monographs (see e.g., \cite{Rankin1, Schoeneberg, Koblitz}).

Denote by $\GN$ the discrete subgroup of $\text{SL}_2(\mathbb{R})$
and by $\UH$ the upper half-plane. We call a meromorphic function $f$ on 
$\mathcal{H}$ a {\it meromorphic modular form} of {\it weight k} for 
$\Gamma$ if
\begin{equation*}
f\Big(\frac{az+b}{cz+d}\Big)=(cz+d)^kf(z)\,, \qquad 
\begin{pmatrix} a & b \\ c & d\end{pmatrix} \in \Gamma\,,
\end{equation*}
for all $z\in\mathcal{H}$, and 
$f$ is meromorphic at the cusps (i.e., $\mathbb{R} \cup \{\infty\}.$)
If $k=0,$ then $f$ is called a {\it modular function} on $\GN.$ Further, we 
say that $f$ is a 
{\it holomorphic modular form} if $f$ is holomorphic on $\mathcal{H}$ and 
holomorphic at the cusps. A holomorphic modular form is said to be a 
{\it cusp form} if it vanishes at each cusp of $\GN.$ 

As is customary, we denote by $\Mk$ (resp. $\Sk$) the space of holomorphic 
modular forms (resp. cusp forms) of weight $k$ for $\GN.$ Moreover, we denote 
by $\Minf$ the space of weakly holomorphic modular forms on $\GN$ 
(i.e. holomorphic on $\UH$ but not necessarily at the cusps). For example, 
$E_4(z)\in M_4(\SL)$ and $E_6(z) \in M_6(\SL)$.
It is well known that $E_2(z)$ is not a modular form of weight $2$ 
but is often referred to as a {\it quasi}-modular form
of weight 2. In fact, there is no modular form of weight $2$  
for $\SL.$ Also recall that the two well known modular forms for $\SL$ namely, 
the {\it Delta} function and the {\it modular} $j$-function
are defined as 
\begin{equation*}
\D(z):=\frac{E_4^3(z)-E_6^2(z)}{1728}= \, q\prod_{n=1}^\infty(1-q^n)^{24} 
:= \, \sum_{n=1}^\infty\tau(n)q^n\,, \qquad \tau(n) \in \mathbb{Z}\,,
\end{equation*} 
\begin{equation*}
j(z):=\frac{E_4^3(z)}{\D(z)}=\frac{1}{q}+744+196884q+21493760q^2+\cdots.
\end{equation*}
Since the only cusp for $\SL$ is at infinity (i.e., $q=0$), then clearly 
$\D\in S_{12}(\SL)$ and $j\in M_0^{\infty}(\SL)$. In fact, all meromorphic
modular functions on $\SL$ are rational functions of $j$ (see \cite{Rankin1}, \cite{Ono}).
The integer coefficients $\tau(n)$ above, are the famous tau-functions of Ramanujan, 
who proved and conjectured many of its properties. The characterization of $\tau(n)$ is 
one of the important problems in number theory. For instance, the well-known
Ramanujan conjecture: $|\tau(n)| < n^{11/2}\sigma_0(n)$ (where $\sigma_0(n)$ is
the number of divisors of $n$) was proved by Deligne in 1974. In the present context,
it is also relevant to point out the relationship between $\D(z)$ and the particular
modular solution of the Chazy equation \eqref{chazy}, namely
\begin{equation}
y(z) = \pi i E_2(z) = \frac{1}{2}\frac{\D'(z)}{\D(z)}
\label{D-chazy}
\end{equation}
Equation \eqref{D-chazy} follows by taking the logarithmic derivative of
the infinite product formula for $\D(z)$ above, and comparing it with the
$q$-expansion for $E_2(z)$ in \eqref{E}.

For the congruent subgroup $\GN=\Gamma_0(2)$ of $\SL$ there are two 
cusps, namely zero and infinity. In this case, it is known that for 
$k$ even and $k \geq 0, \, \dim (M_k(\Gamma_0(2)) = [k/4]+1$
while $\dim (S_k(\Gamma_0(2)) = [k/4]-1$ (see e.g. \cite{Rankin1}),
here $[n]$ denotes the integer part of $n$.
Thus unlike the space $M_2(\SL)$, which is trivial, the space of modular 
forms of weight $2$ on $\Gamma_0(2)$ is 1-dimensional. In fact, it was
recently proved that $\mEt_2(z)=\mPt(q)$ is the ``unique" modular form of 
weight $2$ on $\Gamma_0(2)$ \cite[Lemma 3.3]{BBY}. Yet another recent result 
is that for $k\geq 4,$ the functions $\mE_k(z)$ defined in \eqref{eisen},
are the modular forms of weight $k$ on $\Gamma_0(2)$ which vanish at the cusp 
zero \cite[Theorem 1.1]{Boylan}. 
When $k=2$, it turns out that $\mE_2(z)=\mP(q)$ is a quasi-modular 
form on $\Gamma_0(2)$ (see \eqref{rel-mE2} in \S \ref{sym}), like $E_2(z)$ is for
for the modular group $\SL$. $\mE_2(z)$  plays an important role in the theory of 
modular forms of level 2, as well as in our subsequent discussions.

Define a modular form of weight $4$ on $\Gamma_0(2)$ by
\begin{equation}\label{D}
\mD(z):=\frac{\mEt_2^2(z)-\mE_4(z)}{64}=q+8q^2+28q^3+64q^4+ \cdots.
\end{equation} 
It is easy to check that $\mD$ vanishes at the cusp infinity, but not at 
the cusp zero ($q=1$). Hence $\mD$ is not a cusp form, but $\mD(z)\mE_4(z)$ is 
a cusp form of weight 8. (Since $\dim (S_k(\Gamma_0(2)) = [k/4]-1$, the smallest weight 
of a cusp form on $\Gamma_0(2)$ is 8.) 


Analogous to the Ramanujan tau-functions, the coefficients of $\mD$ have number theoretic
significance. They are related to the number of representations 
of integers as sums of triangular numbers $t_n := n(n+1)/2,\, n=0, 1, 2, \dots$.
By using Ramanujan's theory of elliptic functions, Hahn \cite[Theorem 2.4]{Hahn-divisor}
has proved that
\begin{equation*}
\mD(z)= q\sum_{n=0}^{\infty}\delta_8(n)q^n \,,
\end{equation*}
where $\delta_m(n)$ the number of representations of $n$ as
the sum of $m$ triangular numbers, for $m\geq 1$.
The problem of finding explicit formulas for $\delta_m(n)$ is one of the 
classical problems in combinatorics and number theory \cite{Williams}.

Next, we define a modular function of weight 0 on $\Gamma_0(2)$ by  
\begin{equation}
j_2(z):=\frac{\mEt_2^2(z)}{\mD(z)}=
\frac{1}{q}+40+276q-2048q^2+\cdots \quad \in M_0^{\infty}(\Gamma_0(2))\,,
\label{j2}
\end{equation}
which generates
the field of modular functions on $\Gamma_0(2)$. 
We will show in this paper that in a certain sense the functions $\mD$ 
and $j_2$ play similar roles for the congruence subgroup
$\Gamma_0(2)$ as is done by $\Delta$ and $j$ respectively, for $\SL$.
However, an important difference is that $\Delta$ is a cusp form for $\SL$, whereas
$\mD$ is not a cusp form for $\Gamma_0(2)$. 


\section{Symmetry transformation}\label{sym}
It is well known that there exist both algebraic as well as differential
identities among the Eisenstein series \eqref{Ek} associated with the modular
group $\SL$ (see e.g., \cite[pp 123]{Koblitz}). These identities can be derived 
essentially from the fact that the vector space $M_k(\SL)$ is 
finite dimensional, and that $M_kM_l \subset M_{k+l}$.
The latter fact is a direct consequence 
of the transformation property of modular forms. In fact, the differential identities 
\eqref{dPQR} of Ramanujan can be proved by re-interpreting them as the known 
identities among the Eisenstein series $E_2(z),\, E_4(z),\, E_6(z)$, and $E_8(z)$.
In this section, we demonstrate that there also exist similar 
identities among the Eisenstein series \eqref{eisen} and the function $\mEt_2(z)$
associated with the subgroup $\GN_0(2)$ of the modular group. Furthermore,
Ramamani's differential system \eqref{dmpq} can be directly deduced from these identities.

Note that $\mEt_2(z)$ is {\em not} an Eisenstein series for $\GN_0(2)$ but it can
be expressed as a linear combination of Eisenstein series (see \eqref{EEE} below). 
To establish this result, it is convenient to first introduce an alternative representation for 
$\mEt_2(z)$. By the elementary fact
$$\frac{x}{1+x}=\frac{x}{1-x}-\frac{2x^2}{1-x^2}\,,$$ 
we find from \eqref{Et} and the definition of $\mPt(q)$ in \eqref{mpptq} that
{\allowdisplaybreaks\begin{align}
\mEt_2(z)&=1+24\sum_{n=1}^{\infty}\frac{nq^n}{1+q^n} \nonumber\\
&=1+24\sum_{n=1}^{\infty} \frac{nq^n}{1-q^n}-24\sum_{n=1}^{\infty}\frac{2nq^{2n}}{1-q^{2n}} 
=1+24\sum_{n=1}^{\infty}\frac{(2n-1)q^{2n-1}}{1-q^{2n-1}}
\label{alt-e} \,,
\end{align}}
where we decomposed the first sum in the left hand side of
\eqref{alt-e} into odd and even parts to arrive at the final result. 
\begin{lemma} \label{3mp2mpt}
Let $E_2(z)$, $\mE_2(z)$, and $\mEt_2(z)$ be defined as in \eqref{E}, 
\eqref{eisen}, and \eqref{Et}, then we have that
\begin{equation}\label{EEE}
E_2(z)=3\mE_2(z)-2\mEt_2(z).
\end{equation}
\end{lemma}

\begin{proof}
First, recall that $\mE_2(z) = \mP(q)$. Then, from the definition of $\mPt(q)$ in 
\eqref{mpptq} and by \eqref{alt-e} we obtain
 {\allowdisplaybreaks\begin{align*}
3\mE_2(z)-2\mEt_2(z)=&3\Big(1+8\sum_{n=1}^{\infty}\frac{(2n-1)q^{2n-1}}{1-q^{2n-1}}
-8\sum_{n=1}^{\infty}\frac{2nq^{2n}}{1-q^{2n}}\Big)
-2\Big(1+24\sum_{n=1}^{\infty}\frac{(2n-1)q^{2n-1}}{1-q^{2n-1}}\Big)\\
=&1-24\sum_{n=1}^{\infty}\frac{(2n-1)q^{2n-1}}{1-q^{2n-1}}-24\sum_{n=1}^{\infty}
\frac{2nq^{2n}}{1-q^{2n}}
=1-24\sum_{n=1}^{\infty}\frac{nq^n}{1-q^n}=E_2(z)\,.
\end{align*}}
The last identity for $E_2(z)$ follows after expanding $(1-q^n)^{-1}$ in
a geometric series and then using \eqref{E}. 
\end{proof}
Comparing the above representation for $E_2(z)$ with the left hand 
side of \eqref{alt-e}, yields
\begin{equation*} 
\mEt_2(z)=2E_2(2z)-E_2(z),
\end{equation*}
which we then put in \eqref{EEE} to eliminate $\mEt_2(z)$, and obtain the following
expression 
\begin{equation}\label{EE}
\mE_2(z)=\frac{4}{3}E_2(2z)-\frac{1}{3}E_2(z).
\end{equation}
>From \eqref{EE} and the transformation for $E_2(z)$ which is a quasi-modular form
of weight 2 on $\SL$, we can derive the transformation formula for $\mE_2(z)$.

\begin{lemma} 
For $\left(\begin{smallmatrix} a & b \\ c & d\end{smallmatrix}\right)\in \Gamma_0(2),$ 
$\mE_2(z)$ transforms like a quasi-modular form of weight 2 as follows
\begin{equation}\label{rel-mE2}
\mE_2 \Big( \frac{az+b}{cz+d} \Big)=(cz+d)^2\mE_2(z)
+\frac{2}{\pi i}c(cz+d).
\end{equation}
\end{lemma}
\begin{proof}
Recall for $\left(\begin{smallmatrix} a & b \\ c & d\end{smallmatrix}\right)\in \SL,$ 
the transformation formula \cite[pp 68]{Schoeneberg} 
\begin{equation}\label{rel-E2}
E_2\Big( \frac{az+b}{cz+d} \Big) =(cz+d)^2E_2(z)+\frac{6}{\pi i}c(cz+d).
\end{equation}
Then we have for $\left(\begin{smallmatrix} a & b \\ c & d\end{smallmatrix}\right)\in 
\Gamma_0(2),$ 
$$E_2\left(2\Big(\frac{az+b}{cz+d}\Big)\right)=
E_2 \left(\frac{a(2z)+2b}{\frac{c}{2}(2z)+d}\right)=(cz+d)^2E_2(2z)
+\frac{3}{\pi i}c(cz+d).$$ Hence, by \eqref{EE}, we have that
{\allowdisplaybreaks\begin{align*}
\mE_2 \Big( \frac{az+b}{cz+d} \Big)=&\frac{4}{3}\Big((cz+d)^2E_2(2z)+
\frac{3}{\pi i}c(cz+d)\Big)
-\frac{1}{3}\Big((cz+d)^2E_2(z)+\frac{6}{\pi i}c(cz+d)\Big)\\
=&(cz+d)^2\mE_2(z)+\frac{2}{\pi i}c(cz+d).
\end{align*}}
\end{proof}
The above lemma allows us to establish the following
transformation property for modular forms of weight $k$ on $\Gamma_0(2)$.
\begin{lemma}
Let $f \in M_k(\Gamma_0(2))$, then 
$\delta f-\displaystyle\frac{k}{4}\mE_2f \in M_{k+2}(\Gamma_0(2))$
where the differential operator 
$\delta := \displaystyle q\frac{d}{dq} = \frac{1}{2\pi i}\frac{d}{dz}$.
\label{transf} \end{lemma}
\begin{proof}
For $\left(\begin{smallmatrix} a & b \\ c & d\end{smallmatrix}\right)
\in \Gamma_0(2),$ we have  
$$f'\left(\frac{az+b}{cz+d}\right)=
(cz+d)^2\left((cz+d)^kf(z)\right)'=(cz+d)^{k+2}f'(z)+kc(cz+d)^{k+1}f(z)\,,$$
and from the transformation \eqref{rel-mE2},
$$ f\Big( \frac{az+b}{cz+d} \Big)\mE_2 \Big( \frac{az+b}{cz+d} \Big)=
(cz+d)^{k+2}f(z)\mE_2(z) +\frac{2}{\pi i}c(cz+d)^{k+1}f(z)\,. $$
Then, by combining the two expressions above yields the desired transformation 
property for $\delta f-\frac{k}{4}\mE_2f$. 
\end{proof}
The above result will be used next to derive the differential identities for
$\mP(q)=\mE_2(z)$, $\mPt(q)=\mEt_2(z)$, and $\mQ(q)=\mE_4(z)$, given by \eqref{dmpq}.
Note that Lemma \ref{transf} does not apply to $\mE_2(z)$ as it is
not a modular form of weight 2. However, a similar calculation as in
Lemma \ref{transf} shows that $\delta \mE_2 - \mE_2^2/4$ is a modular 
form of weight 4. Moreover, since $\dim M_4(\Gamma_0(2)) = 2$, one can write
$\delta \mE_2 - \mE_2^2/4 = a\mEt_2^2 + b\mE_4$ for constant coefficients
$a, b$. By comparing the first 2 terms in the $q$-expansions of both sides
one easily verifies that $a=0$ and $b=-1/4$. Thus we recover the first equation
in \eqref{dmpq} in the form
\begin{equation}
\delta \mE_2 - \frac{\mE_2^2}{4} = -\frac{\mE_4}{4} \,.
\label{dmE2} 
\end{equation}
Applying Lemma \ref{transf} to $\mEt_2$ and $\mE_4$ yields holomorphic
modular forms of weight 4 and 6, respectively on $\Gamma_0(2)$. Then
similar dimensional arguments as above, and comparing the first
few terms of the $q$-expansion on both sides, lead to the remaining equations
in \eqref{dmpq} as follows
\begin{align}
\delta \mEt_2 - \frac{\mE_2\mEt_2}{2} &= -\frac{\mE_4}{2} \,, \label{dmEt2}\\
\delta \mE_4 - \mE_2\mE_4 &= - \mEt_2\mE_4\,. \label{dmE4}
\end{align} 
>From \eqref{D}, and using \eqref{dmEt2}, \eqref{dmE4}, we obtain the relation  
\begin{equation}
\delta \mD = \mE_2\mD \,,
\label{dmD}
\end{equation}
which will be useful in the next section. In fact, \eqref{dmD} plays analogous role for
$\Gamma_0(2)$ (see \eqref{y(z)} below in section 4), as \eqref{D-chazy} does for the 
Chazy equation in the context of the modular group $\SL$. Before closing this section, we 
note the following identity between the $\D(z)$ and $\mD(z)$,
\begin{equation}
\mD^3(z) = \frac{\D^2(2z)}{\D(z)} \,,
\label{DmD}
\end{equation}
which is an easy consequence of re-expressing \eqref{EEE} in terms of the logarithmic
derivatives from \eqref{D-chazy} and \eqref{dmD}, then integrating the resulting equation.
The constant of integration (which is unity) is determined by comparing 
only the first terms of the $q$-expansions for $\D(z)$ and $\mD(z)$.

\section{Main results}

In this section we explore the relation between Ramamani's nonlinear
system \eqref{dmpq} which is equivalent to the equations
\eqref{dmE2}--\eqref{dmE4} for the modular forms on $\Gamma_0(2)$,
and Chazy-type third order nonlinear ODEs studied by Bureau.
We also present a first order system of ODEs that are equivalent to
\eqref{dmpq}, and whose general solution can be prescribed
in terms of the modular function $j_2$. This system may be regarded as
the analogue of the Halphen system \eqref{halphen} for Ramamani's nonlinear
differential system.

In 1987, Bureau \cite{Bureau} investigated a certain class of third order nonlinear
ODEs, and expressed their general solutions in terms of the Schwarz 
triangle function $s := S(\alpha,\beta,\gamma;z)$ which satisfies
\begin{equation}
\frac{s'''}{s'}-\frac{3}{2}\left(\frac{s''}{s'}\right)^2 + 
\frac{s'^2}{2} V(s) = 0\,,
\label{schwarz}
\end{equation}
and where $V(s)$ is given by
\begin{equation}
V(s)=\frac{1-\alpha^2}{s^2}+\frac{1-\beta^2}{(s-1)^2}+
\frac{\alpha^2+\beta^2-\gamma^2-1}{s(s-1)} \,.
\label{Vdef}
\end{equation}
The Schwarzian equation \eqref{schwarz}
describes the conformal mappings of the upper-half $s$-plane to
the interior of a region in the extended complex plane and
bounded by three regular circular arcs.  If $\alpha$, $\beta$, and $\gamma$ are
non-negative real numbers such that $\alpha+\beta+\gamma<1$, then the
angles subtended at the vertices $s=0$, $s=1$, and $s=\infty$
of this triangle are $\alpha\pi$, $\beta\pi$, and $\gamma\pi$, respectively.
Furthermore, if $\alpha$, $\beta$, and $\gamma$ are chosen to be
either zero or reciprocals of positive integers, then $s(z)$ is an 
invertible, meromorphic function on the interior of a circle in the extended
complex plane, and cannot be analytically continued across this circle, which 
turns out to be a natural barrier (see e.g., \cite{Nehari}).

The solution to \eqref{schwarz} is given implicitly in terms
of the inverse function $z(s)$ which exists when 
$\alpha, \beta, \gamma \in \{0\} \cup \{1/n\, : n \in \mathbb{N}\}$.
It can be shown that (see e.g., \cite{Nehari, Ford}) the inverse function
is given by the ratio $z(s)=u_1(s)/u_2(s)$,
where $u_1(s)$ and $u_2(s)$ are any two independent solutions to the Fuchsian 
differential equation
\begin{equation*}
u''+\frac{V(s)}{4}u=0 \label{ueqn}
\end{equation*}
with regular singular points at $s=0, \, s=1$ and $s = \infty$.
If we set
\begin{equation*}
 \chi(s) = s^{(\alpha-1)/2}(s-1)^{(\beta-1)/2}u(s)\,,
\end{equation*}
then $\chi(s)$ satisfies the hypergeometric equation,
\begin{equation}
s(s-1)\chi''+[(a+b+1)s-c]\chi'+ab\chi=0,
\label{hyper}
\end{equation}
where $a=(1-\alpha-\beta+\gamma)/2$, $b=(1-\alpha-\beta-\gamma)/2$, 
and $c=1-\alpha$. Thus, in terms of any fundamental set of solutions
$\chi_1,\,\chi_2$ of the hypergeometric equation \eqref{hyper}, the 
general solution to \eqref{schwarz} is given by the inverse $s(z)$
of the following function of $s$
\begin{equation}
z(s) = \frac{A\chi_1(s)+B\chi_2(s)}{C\chi_1(s) +D\chi_2(s)}\,,
\label{zdef}
\end{equation}
where $A,B,C,D$ are complex constants with $AD-BC \neq 0$.
Note that a different choice
$\widetilde{\chi}_1=a_1\chi_1+a_2\chi_2$ and $\widetilde{\chi}_2=b_1\chi_1+b_2\chi_2$
for the pair of independent solutions of \eqref{hyper} induces
a linear fractional transformation of $z(s)$, but leaves invariant
the general form of the equation \eqref{zdef}. That is,
\begin{equation}
\widetilde{z}(s) = \frac{A\widetilde{\chi}_1(s)+B\widetilde{\chi}_2(s)}
{C\widetilde{\chi}_1(s) +D\widetilde{\chi}_2(s)}
= \frac{\widetilde{A}\chi_1(s)+\widetilde{B}\chi_2(s)}
{\widetilde{C}\chi_1(s) +\widetilde{D}\chi_2(s)} = \frac{az(s)+b}{cz(s)+d}\,,
\label{ztransf}
\end{equation}
where $a,\,b,\,c,\,d$ are complex constants with $ad-bc \neq 0$.
When the inverse of $z(s)$ exists, it follows from \eqref{zdef} and \eqref{ztransf} 
that if $s(z)$ is a solution
of the Schwarzian equation \eqref{schwarz}, then the function
\begin{equation}
\widetilde{s}(z) := s\Big(\frac{az+b}{cz+d}\Big)\,, \qquad
\begin{pmatrix} a & b \\ c & d\end{pmatrix} \in
\text{SL}_2(\mathbb{C})\,,
\label{stransf}
\end{equation}
is also a solution of \eqref{schwarz}. Hence, starting from a particular solution
namely, the Schwarz triangle function $s(z)=S(\alpha,\,\beta,\,\gamma;\,z)$ 
with $\alpha, \beta, \gamma \in \{0\} \cup \{1/n\, : n \in \mathbb{N}\}$, the general 
solution of \eqref{schwarz} can be constructed via the transformation 
property \eqref{stransf}.

In his work, Bureau considered all third order scalar ODEs whose general
solutions are given in terms of the Schwarz triangle functions $S(\alpha,\beta,\gamma;z)$
where the parameters $\alpha, \beta, \gamma \in \{0\} \cup \{1/n\, : n \in \mathbb{N}\}$.
In particular, Bureau's class of ODEs includes the Chazy equation \eqref{chazy} 
that corresponds to the Schwarz triangle function $S(1/2,1/3,0;z)$ which
is the same as the modular $j$-function for $\SL$ \cite{Ford}. 
In this section, we establish 
similar results for the congruence subgroup $\Gamma_0(2)$. We begin by rewriting
\eqref{dmE2}--\eqref{dmE4} into a single equation in order to deduce the following
result.

\begin{theorem}\label{y}
Let $y(z):= \pi i\mE_2(z)$, then $y(z)$ satisfies the 
third order nonlinear differential equation
\begin{equation}\label{y(z)}
y'''= 2yy''-y'\,^2 + 2\frac{(y''-yy')^2}{2y'-y^2}\,.
\end{equation}
\end{theorem}

\begin{proof}
The proof involves elimination of $\mEt_2$ and $\mE_4$ using 
\eqref{dmE2}--\eqref{dmE4}, and then obtaining a single equation 
for $\mE_2$. We briefly outline the steps of this calculation
below.
 
Solving for $\mEt_2$ from \eqref{dmE4}, we have 
$\mEt_2 = \mE_2-\delta\mE_4/\mE_4$. Now applying the differential operator
$\delta$ to both sides of this expression, and using \eqref{dmEt2} yields,
$$ \delta\Big(\mE_2-\frac{\delta\mE_4}{\mE_4}\Big) = 
\frac{\mE_2}{2}\Big(\mE_2-\frac{\delta\mE_4}{\mE_4}\Big) -\frac{\mE_4}{2}\,. $$
If we use \eqref{dmE2} to eliminate $\mE_4$ from the
equation above, then it is clear that we obtain a third order ODE for $\mE_2$.
Next we put $\mE_2 = y/\pi i$ into this third order ODE, and after
straightforward calculations, arrive at the desired result.
\end{proof}

Rankin in \cite{Rankin}, showed that $\Delta(z)$ 
associated with the modular group $\SL$, satisfies an ODE which is
homogeneous of degree $4$ (in both $\Delta$ and its derivatives), namely 
$$
2\D''''\D^3-10\D'''\D'\D^2-3\D''\,^2\D^2+
24\D''\D'\,^2\D-13\D'\,^4=0\,. $$
The above equation also follows directly from the Chazy equation 
\eqref{chazy}, if we express $y$ as the logarithmic derivative of $\Delta$
as in \eqref{D-chazy}. We prove a similar result for $\GN_0(2)$.

\begin{corollary}\label{coro}
The modular form $\mD(z)$ defined as in \eqref{D}, satisfies the
following ODE, which is of degree 6 in both $\mD$ and its derivatives.
{\allowdisplaybreaks\begin{align}
&\mD''''(8\mD''\mD^4-10\mD'\,^2\mD^3)+8\mD'''\,^2\mD^4 + \mD'''(10\mD'\,^3\mD^2+16\mD''\mD'\mD^3)
\nonumber \\ 
&-20\mD''\,^3\mD^3-48\mD''\,^2\mD'\,^2\mD^2 -60\mD''\mD'\,^4\mD+25\mD'\,^6 =0
\label{Dode}  
\end{align}}
\end{corollary}
\begin{proof}
>From \eqref{dmD}, the solution $y(z)=\pi i \mE_2(z)$ in Theorem \ref{y}
can be rewritten as $y = \mD'/2\mD$. Substituting this expression for $y$ 
in equation \eqref{y(z)}, we obtain the desired result after a lengthy but 
straightforward calculation.
\end{proof}


We point out that \eqref{y(z)} also appears in Bureau's paper
\cite[pp 352, Eq (55) with n=1/2, m=p=0]{Bureau}. These parameter
values correspond to those of the underlying Schwarz triangle function
which in fact, is a modular function on $\Gamma_0(2)$. Indeed, if we define
the $\Gamma_0(2)$-invariant function
\begin{equation}
s(z):=\frac{\mEt_2^2(z)}{64 \mD(z)}=\frac{j_2(z)}{64}\,,
\label{s(z)}
\end{equation}
then by using equations \eqref{dmE2}--\eqref{dmD}
and \eqref{D}, it is possible to express $\mEt_2,\,\mD,\,\mE_4$ 
and $\mE_2$ in terms of $s(z)$ and its derivatives, as given below.
\begin{equation}
\mEt_2=\frac{\delta s}{1-s}\,, \quad  
\mD=\frac{(\delta s)^2}{64s(s-1)^2}\,, \quad  
\mE_4=\frac{(\delta s)^2}{s(s-1)}\,, \quad  
\mE_2 = \delta \ln\Big(\frac{(\delta s)^2}{s(s-1)^2}\Big)\,. \quad  
\label{mEs}
\end{equation}
Using the expressions for $\mE_4$ and $\mE_2$ from \eqref{mEs} in \eqref{dmE2}
yields the following result.  
\begin{lemma} \label{s}
The $\Gamma_0(2)$-invariant modular function $s(z)$ defined 
as in \eqref{s(z)} satisfies \eqref{schwarz} with parameter values
$\alpha=1/2, \, \beta=0, \, \gamma = 0$ in \eqref{Vdef}. 
\end{lemma}
We emphasize that the modular function $s(z)$ in \eqref{s(z)} is only a special
solution of \eqref{schwarz} with the parameter values given in Lemma \ref{s}.
The general solution can be obtained via the linearization procedure 
leading to \eqref{zdef}, as described earlier.
In this case, the parameter values of the associated hypergeometric 
equation \eqref{hyper} are $a=b=1/4$ and $c=1/2$. Thus for example, 
if in \eqref{zdef}, we choose 
$$
\chi_1(s)=\,_2F_1(\frac{1}{4},\frac{1}{4},\frac{1}{2};s), \qquad 
\chi_2(s)= \sqrt{s}\,\,_2F_1(\frac{3}{4},\frac{3}{4},\frac{3}{2};s), \quad |s| < 1\,,
$$ 
as the two independent hypergeometric solutions, then we obtain (after appropriate
analytic continuation of $\chi_1(s)$ and $\chi_2(s)$ in the $s$-plane and inverting
$z(s)$) the general solution $s(z)$ of \eqref{schwarz} with
$\alpha=1/2, \beta=0, \gamma = 0$. Here $_2F_1(a,b,c;s)$ denotes the standard
hypergeometric series. From the general solution $s(z)$, one can also
obtain the general solution to \eqref{y(z)} in Theorem \ref{y}. 
The key idea is that if instead of $S(1/2,0,0;z)$, the general
solution $s(z)$ is used to express $\mEt_2,\,\mD,\,\mE_4$
and $\mE_2$ in \eqref{mEs}, then
the resulting complex functions $\mEt_{2}^{C},\,\mD^C,\,\mE_4^C$
and $\mE_2^C$ will still satisfy the differential relations 
\eqref{dmE2}--\eqref{dmD} even though these functions are no longer
modular forms on $\Gamma_0(2)$ (except when $s=S(1/2,0,0;z)$).

\begin{theorem} \label{ygen}
Let $s(z)$ be the general solution of 
\eqref{schwarz} with $\alpha=1/2, \beta=0, \gamma = 0$, then
\begin{equation}
y(z) = \frac{1}{2}\Big[ \frac{s''}{s'}-\Big(\frac{1}{2s}+
\frac{1}{s-1}\Big)s'\Big]  
\label{yg}
\end{equation}
is the general solution to \eqref{y(z)}. Moreover, \eqref{y(z)} has
the following transformation property:\, if $y(z)$ is a solution then
\begin{equation}
\widetilde{y}(z) = \frac{1}{(cz+d)^2}y\Big(\frac{az+b}{cz+d}\Big)
-\frac{c}{cz+d}\,, \qquad 
\begin{pmatrix} a & b \\ c & d\end{pmatrix} \in \text{SL}_2(\mathbb{C}).
\label{ytransf}
\end{equation}
is also a solution of \eqref{y(z)}.
\end{theorem}
 
\begin{proof}
We define the function $\mE_2^C$ in \eqref{mEs} using the general 
solution $s(z)$, and set 
$y(z) = (\pi i/2)\mE_2^C$. Then this gives \eqref{yg}. Now, repeating 
the same calculations as in 
Theorem \ref{y} by using \eqref{dmE2}--\eqref{dmE4} with the functions 
$\mEt_{2}^{C},\,\mE_4^C \, \mE_2^C$, leads to the first part of the Lemma. 
Alternatively, one can directly verify the assertion by substituting 
the expression for $y(z)$ from \eqref{yg} into \eqref{y(z)}, and making use
the Schwarzian equation \eqref{schwarz}.

The transformation property \eqref{ytransf} follows immediately from the
transformation of $s(z)$ given by \eqref{stransf}, if we use $\widetilde{s}(z)$
from \eqref{stransf} to define the function $\widetilde{y}(z)$ in \eqref{yg}.
\end{proof}

In 1881, Halphen also considered a more general version of the nonlinear differential
system \eqref{halphen}, and presented its solution in terms of hypergeometric functions 
\cite{Halphen1}. Recently, in \cite{Ablowitz-H1,Ablowitz-H2}, the authors 
rediscovered this
system arising as a special case of self-dual Yang-Mills equations in mathematical physics, 
and called it the generalized Darboux-Halphen (gDH) system. This system is given by
\begin{align} 
u_1'&=u_2u_3-u_1(u_2+u_3)+\tau^2\,,  \nonumber\\ 
u_2'&=u_3u_1-u_2(u_3+u_1)+\tau^2\,,  \label{ghalphen}\\ 
u_3'&=u_1u_2-u_3(u_1+u_2)+\tau^2\,,  \nonumber  
\end{align}
$$ \tau^2 = \alpha^2(u_1-u_2)(u_2-u_3)+\beta^2(u_2-u_1)(u_1-u_3)+\gamma^2(u_3-u_1)(u_2-u_3) \,,$$
for functions $u_i(z) \neq u_j(z),\, i\neq j, \, i,j=1,2,3$, and complex constants 
$\alpha,\, \beta,\,\gamma$. The system \eqref{ghalphen} reduces to
\eqref{halphen} when $\alpha = \beta = \gamma = 0$.
Note also that the gDH system \eqref{ghalphen} appears in a slightly
different form from the original system considered by Halphen 
\cite[pp 1405, Eq (5)]{Halphen1}, but it can be transformed into Halphen's 
system by suitable linear combinations of the variables $u_1(z),\, u_2(z), \, u_3(z)$.

In \cite{Ablowitz-H1,Ablowitz-H2}, the solution to the gDH 
system was obtained by parametrizing the variables $u_i(z)$ in terms of the 
solution (and its derivatives)
of the Schwarzian equation \eqref{schwarz} discussed earlier. In this case, the complex constants
$\alpha,\, \beta,\,\gamma$ in \eqref{Vdef} are the same as those in the function
$\tau^2$ above. We now derive a gDH system for the congruence subgroup $\Gamma_0(2)$,
and show its equivalence to Ramamani's system \eqref{dmpq}. 

With the aid of \eqref{schwarz} and \eqref{stransf}, one can
verify the following statements.
\begin{lemma}\label{dh}
Let $u_1(z),\, u_2(z), \, u_3(z)$ be defined by
\begin{equation}
u_1 = -\frac{1}{2}\Big[\ln\Big( \frac{s'}{s} \Big)\Big]'\,, \quad
u_2 = -\frac{1}{2}\Big[\ln\Big( \frac{s'}{s-1} \Big)\Big]'\,, \quad
u_3 = -\frac{1}{2}\Big[\ln\Big( \frac{s'}{s(s-1)} \Big)\Big]'\,,
\label{udef}
\end{equation}
where $s(z)$ is the $\Gamma_0(2)$-invariant modular function in \eqref{s(z)}, then
\eqref{udef} solves the gDH system \eqref{ghalphen} with $\alpha = 1/2$ and $\beta= \gamma=0$.
Moreover, if $u_i(z), \, i=1,2,3$ are solutions of \eqref{ghalphen} then so are 
\begin{equation*}
\widetilde{u_i}(z) = \frac{1}{(cz+d)^2}u_i\Big(\frac{az+b}{cz+d}\Big)
+\frac{c}{cz+d}\,,   \qquad
\begin{pmatrix} a & b \\ c & d\end{pmatrix} \in \text{SL}_2(\mathbb{C}).
\end{equation*}
\end{lemma}
The $\Gamma_0(2)$-invariance of $s(z)$ implies that
\begin{equation*}
s\Big(\frac{az+b}{cz+d}\Big) = s(z), \quad \text{and} \quad  
s'\Big(\frac{az+b}{cz+d}\Big) = (cz+d)^2s'(z)\,, \qquad 
\begin{pmatrix} a & b \\ c & d\end{pmatrix} \in \Gamma_0(2)\,. 
\end{equation*}
Then it follows from \eqref{udef} that the gDH variables $u_i(z)$ are quasi-modular 
forms of weight 2 on $\Gamma_0(2)$, and they transform according to 
\begin{equation*}
u_i\Big(\frac{az+b}{cz+d}\Big) = (cz+d)^2u_i(z) - c(cz+d)\,, \qquad 
\begin{pmatrix} a & b \\ c & d\end{pmatrix} \in \Gamma_0(2)\,,
\end{equation*}
for $i=1,2,3$. Note however that the differences $u_i - u_j\,, \,\, i\neq j$ are forms of weight 2 
on $\Gamma_0(2)$, and that the modular function $s(z)$ itself can be expressed by the
cross-ratio
\begin{equation*}
s(z) = \frac{u_1-u_3}{u_1-u_2}\,.
\end{equation*}
It is worth noting here that the variables $u_i$ associated with the DH system \eqref{halphen} 
can also be represented similarly as in Lemma \ref{dh} in terms of 
\begin{equation}
\lambda(z) := \frac{\vartheta_2^4(0|z)}{\vartheta_3^4(0|z)}\,,
\label{lambda}
\end{equation}
which is the modular function for the subgroup 
$\GN(2) := \{\gamma \in \SL:\gamma \equiv I \, \pmod{2}\}$ \cite{Takhtajan}.
In fact, by replacing $s$ by $\lambda$ in \eqref{udef} and using certain identities 
among the null theta functions (see below), one can recover Halphen's solution for \eqref{halphen}.

>From \eqref{mEs} and \eqref{udef}, it is also
possible to express the modular forms $\mEt_2, \mE_2$ and $\mE_4$ in terms
of the DH variables $u_i$ as follows
\begin{equation*}
\pi i \mEt_2 = u_1-u_3, \quad \pi i \mE_2 = -(u_2+u_3), \quad 
\pi^2\mE_4 = (u_1-u_3)(u_3-u_2)\,.
\end{equation*}
Conversely, $u_i(z),\, i=1,2,3$ are rational functions of $\mEt_2, \mE_2, \mE_4$. Thus there
exists a bijection between the solutions of the gDH system \eqref{ghalphen} and those 
of \eqref{dmE2}--\eqref{dmE4}, or equivalently, Ramamani's nonlinear system \eqref{dPQR}.

There is also an alternative representation of the solutions \eqref{udef} of
the gDH system \eqref{ghalphen} in terms of null theta functions, like in the
case of the DH system \eqref{halphen}. This can be found by first expressing
$s(z)$ as a rational function of the null theta functions, and then
using \eqref{udef} in Lemma \ref{dh} to obtain the resulting expressions for $u_i(z)$.
If we replace $E_2(z)$ and $\mE_2(z)$ in \eqref{EEE} by their equivalent expressions from
\eqref{D-chazy} and \eqref{dmD} respectively, we obtain
\begin{equation*}
 2 \pi i \mEt_2(z) = \frac{3}{2}\frac{\mD'(z)}{\mD(z)} - \frac{1}{2}\frac{\Delta'(z)}{\Delta(z)}
= 2 \frac{\Delta'(2z)}{\Delta(2z)} - \frac{\Delta'(z)}{\Delta(z)}\,,
\end{equation*}
where we have also used \eqref{DmD} to derive the last equality. Equating the expression
for $\mEt_2(z)$ from the first equation in \eqref{mEs} with that in above, and integrating
both sides, yield
\begin{equation}
s(z) = 1 + \frac{1}{64}\frac{\Delta(z)}{\Delta(2z)} \,.
\label{sDelta}
\end{equation}
The constant of integration (1/64) is fixed by considering the $q$-expansion of $\Delta$
and that of $s=j_2/64$ from \eqref{j2}, and comparing the coefficient of $q^{-1}$ on
both sides of \eqref{sDelta}. Next, by using the well known identity $\Delta = \eta^{24}$
where $\eta(z)$ is the Dedekind eta-function, and the following identities 
(see e.g., \cite{WW, Rankin1})
\begin{equation*}
\Delta^{1/8}(z) = \eta^3(z) = \frac{\vartheta_2(z) \vartheta_3(z) \vartheta_4(z)}{2}\,, \qquad
\frac{\eta^2(2z)}{\eta(z)} = \frac{\vartheta_2(z)}{2}\,, \qquad
\vartheta_3^4(z) = \vartheta_2^4(z) + \vartheta_4^4(z)\,,
\end{equation*}
among $\eta(z)$ and the null theta functions $\vartheta_i(0|z):= \vartheta_i(z)$, 
we deduce from \eqref{sDelta} that
\begin{equation}
s(z) = 1 + \frac{1}{64}\left( \frac{\eta(z)}{\eta(2z)} \right)^{24}
= \left( \frac{\vartheta_3^4(z)+\vartheta_4^4(z)}{\vartheta_2^4(z)} \right)^2 \,.
\label{stheta}
\end{equation}
Finally, substituting $s(z)$ from \eqref{stheta} in \eqref{udef}, and making use
of the differential identities
\begin{equation*}
\frac{\vartheta_2'(z)}{\vartheta_2(z)}-\frac{\vartheta_3'(z)}{\vartheta_3(z)}= 
\frac{\pi i}{4}\vartheta_4^4(z)\,, \qquad
\frac{\vartheta_3'(z)}{\vartheta_3(z)}-\frac{\vartheta_4'(z)}{\vartheta_4(z)}= 
\frac{\pi i}{4}\vartheta_2^4(z)\,, \qquad
\frac{\vartheta_2'(z)}{\vartheta_2(z)}-\frac{\vartheta_4'(z)}{\vartheta_4(z)}= 
\frac{\pi i}{4}\vartheta_3^4(z)\,, 
\end{equation*}
we obtain the following expressions for the gDH variables $u_i$ associated with $\GN_0(2)$
\begin{equation}
u_1 = -\frac{1}{2}\Big[\ln\Big( \frac{\vartheta_3^4\vartheta_4^4}
{\vartheta_3^4+\vartheta_4^4} \Big)\Big]', \qquad 
u_2 = -\frac{1}{2}\Big[\ln\Big( \vartheta_3^4+\vartheta_4^4 \Big)\Big]', \qquad
u_3 = -\frac{1}{2}\Big[\ln\Big( \frac{\vartheta_2^8}
{\vartheta_3^4+\vartheta_4^4} \Big)\Big]'\,. 
\end{equation}
Thus, the gDH system \eqref{ghalphen} can be solved in terms of null theta
functions as well. In addition, we also note that there exists a simple relation between 
the modular functions $\lambda(z)$ and $s(z)$ associated with the DH and gDH systems
respectively, namely
\begin{equation*}
s(z) = \left( \frac{2}{\lambda}-1 \right)^2 \,,
\end{equation*}
which follows from \eqref{stheta} and from the definition of $\lambda(z)$ in 
\eqref{lambda}.

\section*{Acknowledgments}
MJA is supported by NSF grant No. DMS 0602151, the work of SC is
supported by NSF grant No. DMS 0307181, and HH is supported in part by
a NSF FRG grant No. DMS 0244660.

\end{document}